\documentclass{amsart}
\usepackage[top=4cm, bottom=4cm, left=3cm, right=3cm]{geometry}

\author{Efe Gürel}
\address{TÜBİTAK Natural Sciences High School, Kocaeli, 41400, Turkey}
\email{efegurel54@gmail.com}

\title[Resolution Of Multiplicative Anomaly]{Resolution Of Multiplicative Anomaly Of Zeta Regularization For Polynomials}

\newtheorem{theorem}{Theorem}[section]

\newtheorem{corollary}[theorem]{Corollary}

\newtheorem{lemma}[theorem]{Lemma}

\DeclareMathOperator*\RegularProd{%
	\mathchoice
	{\ooalign{$\displaystyle\prod$\cr\hidewidth$\displaystyle\coprod$\hidewidth\cr}}%
	{\ooalign{$\textstyle\prod$\cr\hidewidth$\textstyle\coprod$\hidewidth\cr}}%
	{\ooalign{$\scriptstyle\prod$\cr\hidewidth$\scriptstyle\coprod$\hidewidth\cr}}%
	{\ooalign{$\scriptscriptstyle\prod$\cr\hidewidth$%
			\scriptscriptstyle\coprod$\hidewidth\cr}}}

\DeclareMathOperator{\res}{Res}
\DeclareMathOperator{\fp}{\mathcal{F P}}
\DeclareMathOperator{\re}{\mathfrak{Re}}

\usepackage{mathtools}
\usepackage{amsmath}
\usepackage{amssymb}
\usepackage{amsthm}
\usepackage{latexsym}
\usepackage[square,numbers]{natbib}



\begin{document}
\begin{abstract}
    In this paper, the problem of multiplicative anomaly of zeta regularization is solved for polynomials. For a regularizable sequence $\Lambda$, we explicitly calculate the zeta regularized product of $(\Lambda-z_1)\dots(\Lambda-z_n)$ for $z_1,\dots,z_n\in\mathbb{C}$. We give an explicit formula for the discrepancies between polynomials. Our results imply Mizuno's theorem as a special case. We furthermore give novel regularized product formulas for multi-dimensional products in terms of Barnes multiple gamma functions, zeros of the Riemann zeta function and zeros of the Bessel function of the first kind.
\end{abstract}
\maketitle

\section{Introduction}
Let $\Lambda=\{\lambda_k\}_{k\in I}$ be a given sequence of non-zero complex numbers such that the order of the sequence
	\begin{align*}
		 \mu=\inf\left\{ s:\sum_{k\in I} |\lambda_k|^{-s}<\infty, \  s\in \mathbb{R} \right\}
	\end{align*}
	exists, is finite and positive. We define the associated zeta function of the sequence $\Lambda$ as
	\begin{align*}
		\zeta_\Lambda(s)=\sum_{k\in I}\frac{1}{\lambda_k^s}, \qquad(\mathfrak{Re}(s)>\mu),
	\end{align*}
	which denotes an holomorphic function for $\mathfrak{Re}(s)>\mu$. Furthermore, we associate the theta function $\theta_\Lambda$ to the sequence $\Lambda$ as 
	\begin{align*}
		  \theta_\Lambda(t)=\sum_{k\in I}e^{-t\lambda_k}, \qquad (t>0).
	\end{align*}
	Zeta and theta functions are connected by the Mellin-transform formula
	\begin{align*}
		\zeta_\Lambda(s)=\frac{1}{\Gamma(s)}\int_{0}^{\infty}\theta_\Lambda(t)t^{s-1}dt
	\end{align*}
	where $\Gamma$ denotes the gamma function. Throughout, we assume $\theta_\Lambda$ has the following asymptotic expansion as $t\to 0^+$,
	\begin{align}\label{ThetaAsymptoteAxiom}
		\theta_\Lambda(t)\sim\sum_{n=0}^{\infty}c_{i_n}t^{i_n}
	\end{align}
	where $i_n$ is a sequence of increasing real numbers such that $i_0<0$, $i_n\to \infty$ as $n\to \infty$ and $c_{i_n}$ are constants. We note that $-i_0=\mu$. For convenience, let $m=[\mu]$. We define the Weierstrass product associated to the sequence $\Lambda$ by
	\begin{align}\label{WeierstrassDef}
		\Delta_\Lambda(z)=\prod_{k\in I}\left( 1-\frac{z}{\lambda_k} \right)\exp\left( \frac{z}{\lambda_k}+\frac{z^2}{2\lambda_k^2}+\dots+\frac{z^m}{m\lambda_k^m} \right)
	\end{align}
	which is the unique entire function of order $\mu$, zeros at $\Lambda$ and satisfying the relation $        \log\Delta_\Lambda(0)=(\log\Delta_\Lambda)'(0)=\dots=(\log\Delta_\Lambda)^{(m)}(0)=0$.
	\newline
	
	We say that such a sequence $\Lambda$ is regularizable if the function $\zeta_\Lambda(s)$ can be meromorphically continued to a region containing $s=0$ as an holomorphic point. For a regularizable sequence $\Lambda$, the zeta regularized product of $\Lambda$, denoted $\RegularProd_{k\in I} \lambda_k$, is defined as
	\begin{align*} 
		\RegularProd_{k\in I}\lambda_k:=\exp(-\partial_s\zeta_\Lambda(0)).
	\end{align*}
	Existence of the asymptotic expansion \eqref{ThetaAsymptoteAxiom} guarantees that $\Lambda$ is regularizable. Furthermore, by the classical Mellin-continuation lemma \cite{Quine,VorosSpectral} the function $\zeta_\Lambda$ can be meromorphically continued to the whole complex plane with only simple poles at $s=-i_n$. The residues can be determined elegantly in terms of $c_{i_n}$. We will use the fact that $\zeta_\Lambda$ has only simple poles at most. The theory of zeta regularization is mainly focused on the regularization function defined as
	\begin{align*}
		D_\Lambda(z)=\RegularProd_{k\in I}(\lambda_k-z),
	\end{align*}
	where $z\notin\Lambda$. The functions $D_\Lambda$ and $\Delta_\Lambda$ behave very similar. It turns out both are entire functions of order $\mu$ and zeros at $\Lambda$. In fact, they are equal up to an elementary factor (see Theorem \ref{FundamentalThm}). Furthermore, it is more natural to consider $D_\Lambda$ rather than $\Delta_\Lambda$ in theory of functions, given its tame behavior under translation and integral transforms \cite{VorosSpectral}. Lerch \cite{Lerch} famously proved that
	\begin{align}\label{ClassicLetch}
		\RegularProd_{k=0}^{\infty}(k+x)=\frac{\sqrt{2\pi}}{\Gamma(x)}
	\end{align}
	which is known as Lerch's formula in the literature. Lerch originally proved  $\partial_s\zeta(0,x)=\log\frac{\Gamma(x)}{\sqrt{2\pi}}$ where $\zeta(s,x)$ is the Hurwitz zeta function, the zeta function associated to the sequence $\mathbb{N}+x$. This formula seems to be the first example of evaluation of a regularization function. Let $\mathcal{FP}f(s)$ be the constant term of the Laurent expansion of the function $f(z)$ near $z=s$ and $H_n$ be the $n$th harmonic number. Lerch's result can be generalized to give the fundamental theorem of zeta regularization.
	\begin{theorem}\label{FundamentalThm}
		If $\Lambda$ is regularizable, then so is $\Lambda-z$ where $z\notin \Lambda$. The function $D_\Lambda(z)$ can be continued to an entire function of order $\mu$ and zeros at $\Lambda$. Furthermore, the following equation holds,
		\begin{align*}
			D_\Lambda(z)=\exp\left(-\zeta'_\Lambda(0) -\sum_{\ell=1}^{m} \frac{\mathcal{F P}\zeta_\Lambda(\ell)}{\ell}z^\ell-\sum_{\ell=2}^{m}\frac{\res \zeta_\Lambda(\ell) H_{\ell-1}}{k}z^\ell\right)\Delta_\Lambda(z).
		\end{align*}
	\end{theorem}
	On the other hand, Lerch also proved the following not so well-known result
	\begin{align*}
		\RegularProd_{k=0}^{\infty}\big((k+x)^2+y^2\big)=\frac{2\pi}{\Gamma(x+iy)\Gamma(x-iy)}
	\end{align*}
	which is an interesting example of the uncommon equality $\RegularProd a_k b_k=\RegularProd a_k \RegularProd b_k$. Mizuno \cite{Mizuno} gave the following formula for zeta regularization of polynomials.
	\begin{theorem} \label{MizunoTheorem}
		For $z_1,\dots,z_j\in \mathbb{C}$, the following equality holds.
		\begin{align*}
			\RegularProd_{k=0}^{\infty}\prod_{j=1}^{n}(k+z_j)=\prod_{i=j}^{n}\frac{\sqrt{2\pi}}{\Gamma(z_j)}=\prod_{j=1}^{n}\RegularProd_{k=0}^{\infty}(k+z_j)
		\end{align*}
	\end{theorem}
	Theorem \ref{MizunoTheorem} implies that for monic polynomials $\varphi_1,\dots,\varphi_n$ that don't vanish in $\mathbb{Z}^+$, the regulation operation commutes with multiplication:
	\begin{align*}
		\RegularProd_{k=0}^{\infty} \prod_{j=1}^{n} \varphi_j(k)=\prod_{j=1}^{n}\RegularProd_{k=0}^{\infty} \varphi_j(k).
	\end{align*}
	Mizuno also considered a two-dimensional (see Section \ref{SectionMultidim}) and a $q$-analogue of equation \eqref{ClassicLetch} and Theorem \ref{MizunoTheorem}. In general, the equation $\RegularProd a_k b_k=\RegularProd a_k \RegularProd b_k$ does not hold. This phenomena is called the multiplicative anomaly of zeta regularization. For regularizable sequences $\Lambda_1,\dots,\Lambda_n$ and their pointwise product $\Lambda=\Lambda_1\dots\Lambda_n$, we define the discrepancy $F=F(\Lambda_1,\dots,\Lambda_n)$ as
	\begin{align*}
		F=\log \frac{\RegularProd_k \prod_{i=j}^{n} \lambda_{k,j}}{\prod_{j=1}^{n}\RegularProd_k \lambda_{k,j}}=\sum_{j=1}^{n}\zeta'_{\Lambda_j}(0)-\zeta'_\Lambda(0).
	\end{align*}
	By Theorem \ref{MizunoTheorem}, we know that $F(\varphi_1(\mathbb{N}),\dots,\varphi_n(\mathbb{N}))=0$. In this paper, we explicitly evaluate the regularization function given by
	\begin{align*}
			D_\Lambda(z_1,\dots,z_n)=\RegularProd_{k} (\lambda_k-z_1)\dots(\lambda_k-z_n)
	\end{align*}
	for a general sequence $\Lambda$. It turns out that the function $D_\Lambda(z_1,\dots,z_n)$ is an holomorphic function of $n$-variables in the whole complex plane which vanishes whenever there exists a $j$ such that $z_j\in \Lambda$. Furthermore, the functions $\prod_{j=1}^{n} \Delta_\Lambda(z_j)$, $\prod_{j=1}^{n} D_\Lambda(z_j)$ and $D_\Lambda(z_1,\dots,z_n)$ are equal up to an elementary factor.
	We give a closed form of the discrepancy $F(\Lambda-z_1,\dots,\Lambda-z_n)$. Thus, the discrepancy $F(\varphi_1(\Lambda),\dots,\varphi_n(\Lambda))$ can be calculated explicitly for monic polynomials $\varphi_1,\dots,\varphi_n$ that don't vanish in $\Lambda$. This concludes the problem of multiplicative anomaly in the case of polynomials. We apply our results to reobtain the results of Mizuno. We give novel regularization formulas for multi-dimensional products, zeros of the Riemann zeta function and Bessel function of the first kind.
	
	\section{Main Results}
	Throughout this section, we denote $z=(z_1,z_2,\dots,z_n)\in \mathbb{C}^n$ and $s\in\mathbb{C}$. We introduce a family of polynomials that perform an important role in our results. For a non-negative integer $\ell$, define the polynomials $P_\ell(s;z)$ as follows
	\begin{align*}
		P_\ell(s;z)=\sum_{\substack{a_1+\dots+a_n=\ell \\
				a_1,\dots,a_n \in \mathbb{N}}} \binom{-s}{a_1}\dots\binom{-s}{a_n}(-z_1)^{a_1}\dots (-z_n)^{a_n}.
	\end{align*}
	Then trivially $P_0(s;z)=1$ and $P_1(s;z)=\left(\sum_{j=1}^n z_j\right)s$. We need the following crucial lemma.
	\begin{lemma}\label{P_lLemma}
		 For $\ell\ge 2$, let $P_\ell(s;z)=p_{0,\ell}(z)+p_{1,\ell}(z)s+p_{2,\ell}(z)s^2+O\left(s^3\right)$. Then we have $P_\ell(0;z)=p_{0,\ell}(z)=0$,
		\begin{align*}
			\partial_s P_\ell(0;z)=p_{1,\ell}(z)=\sum_{j=1}^{n} \frac{z_j^\ell}{\ell}
		\end{align*}
		and
		\begin{align*}
			p_{2,\ell}(z)=\sum_{j=1}^{n}\frac{H_{\ell-1}}{\ell}z_j^\ell+\sum_{1\le j_1< j_2\le n}\sum_{r=1}^{\ell-1}\frac{z_{j_1}^r z_{j_2}^{\ell-r}}{r(\ell-r)}.
		\end{align*}
	\end{lemma}
	\begin{proof}
		For every $n$-tuple such that $a_1+\dots+a_n=\ell$, there exists at least one $j$ such that $a_j\ge 1$. Therefore $s \big\rvert \binom{-s}{a_j}$ and $s \rvert P_\ell(s;z)$ as polynomials which implies $P_\ell(0;z)=0$. Now we compute the coefficient of $s$ in $P_\ell(s;z)$. If there exists two distinct $a_{j_1},a_{j_2}\ge 1$, then $s^2 \big\rvert \binom{-s}{a_{j_1}}\binom{-s}{a_{j_2}}$ and the $n$-tuple $(a_1,\dots,a_n)$ does not contribute to the coefficient of $s$. Hence we only consider $n$-tuples that have only one non-zero coordinate. All such $n$-tuples are characterized as $a_j=\ell$ for some $j$ and $a_{j'}=0$ for every $j'\neq j$. Since $\binom{-s}{\ell}=\frac{(-1)^\ell}{\ell}s+O\left(s^2\right)$, such an $n$-tuple adds $\frac{z_j^\ell}{\ell}$ to the coefficient of $s$. Putting this all together, we obtain
		\begin{align*}
			\partial_s P_\ell(0;z)=p_{1,\ell}(z)=\sum_{j=1}^{n} \frac{z_j^\ell}{\ell}.
		\end{align*}
		Now we compute $p_{2,\ell}(z)$. If there exists three distinct $a_{j_1},a_{j_2},a_{j_3}\ge 1$ then we have $s^3\big\rvert \binom{-s}{a_{j_1}}\binom{-s}{a_{j_2}}\binom{-s}{a_{j_3}}$ and there is no contribution to $p_{2,\ell}(z)$. Therefore there is at most $2$ non-zero coordinates. We now consider two cases. If there is only one non-zero coefficient, say $a_j=\ell$ for some $j$ and $a_{j'}=0$ for every $j'\neq j$, then this $n$-tuple adds $\frac{H_{\ell-1}}{\ell}z_j^\ell$ since $\binom{-s}{\ell}=\frac{(-1)^\ell}{\ell}s+\frac{(-1)^\ell H_{\ell-1}}{\ell}s^2+O\left(s^3\right)$. If there is two non-zero coefficients, say $a_{j_1}=r$ and $a_{j_2}=\ell-r$ where $a_{j'}=0$ for every $j'\neq j_1,j_2$ and $r=1,\dots,\ell-1$, the $n$-tuple adds $\frac{z_{j_1}^r z_{j_2}^{\ell-r}}{r(\ell-r)}$ since $\binom{-s}{r}\binom{-s}{\ell-r}=\frac{(-1)^\ell}{r(\ell-r)}s^2+O\left(s^3\right)$. Putting all this together we have
		\begin{align*}
			p_{2,\ell}(z)=\sum_{j=1}^{n}\frac{H_{\ell-1}}{\ell}z_j^\ell+\sum_{1\le j_1 < j_2 \le n}\sum_{r=1}^{\ell-1}\frac{z_{j_1}^r z_{j_2}^{\ell-r}}{r(\ell-r)}.
		\end{align*}
		Thus concluding the proof.
	\end{proof}

	\begin{theorem}
		Suppose $\Lambda$ is regularizable. Then $(\Lambda-z_1)\dots(\Lambda-z_n)$ is also regularizable where $z_j\notin \Lambda$ for every $j=1,\dots,n$. 
	\end{theorem}
	
	\begin{proof}
		Define the zeta function associated to the sequence $(\Lambda-z_1)\dots(\Lambda-z_n)$ as
		\begin{align*}
			\zeta_\Lambda(s;{z})=\sum_k \frac{1}{(\lambda_k-z_1)^s \dots (\lambda_k-z_n)^s}.
		\end{align*}
		This series evidently defines an holomorphic function in the region $\re(s)>\frac{\mu}{n}$. Without loss of generality, let $z_1$ have the maximum modulus among $z_1,\dots,z_n$. For $\left|\lambda_k\right|> \left|z_1\right|$, by binomial theorem, we get
		\begin{align*}
			\left(1-\frac{z_j}{\lambda_k}\right)^{-s}=\sum_{\ell=0}^{m} \binom{-s}{\ell} \frac{(-z)^\ell}{\lambda_k^\ell} +O\left(\frac{1}{\lambda_k^{m+1}}\right)
		\end{align*} 
		for every $j=1,\dots,n$. Multiplying the above equation for $j=1,\dots,n$, we obtain
		\begin{align}\label{BinomialMultiplication}
			\prod_{j=1}^{n} \left(1-\frac{z_j}{\lambda_k}\right)^{-s}=\sum_{\ell=0}^{m} \frac{P_\ell(s;{z})}{\lambda_k^\ell} +O\left(\frac{1}{\lambda_k^{m+1}}\right).
		\end{align}
		Therefore we have
		\begin{align}\label{Lambda>ZSeries}
			\begin{split}
				\sum_{\left|\lambda_k\right|> \left|z_1\right|} \frac{1}{(\lambda_k-z_1)^s \dots (\lambda_k-z_n)^s} &-\sum_{\ell=0}^{m}\sum_{\left|\lambda_k\right|> \left|z_1\right|} \frac{P_\ell(s;{z})}{\lambda_k^{ns+\ell}}\\
				&=\sum_{\left|\lambda_k\right|> \left|z_1\right| } \lambda_k^{-ns}\left(\prod_{j=1}^{n} \left(1-\frac{z_j}{\lambda_k}\right)^{-s}-\sum_{\ell=0}^{m} \frac{P_\ell(s;{z})}{\lambda_k^\ell}\right)
			\end{split}
		\end{align}
		where the right hand side of \eqref{Lambda>ZSeries} converges absolutely and uniformly on $\re(s)>\frac{\mu-m-1}{n}$ by equation \eqref{BinomialMultiplication}. By the definition of $\zeta_\Lambda(s)$, we trivially have
		\begin{align}\label{ZetaLambdaComplement}
			\sum_{\ell=0}^{m}\sum_{\left|\lambda_k\right|> \left|z_1\right|} \frac{P_\ell(s;{z})}{\lambda_k^{ns+\ell}}=\sum_{\ell=0}^{m}P_\ell(s;{z})\zeta_\Lambda(ns+\ell)-\sum_{\ell=0}^{m}\sum_{\left|\lambda_k\right|\le \left|z_1\right|} \frac{P_\ell(s;{z})}{\lambda_k^{ns+\ell}}.
		\end{align}
		Decomposing the $\zeta_\Lambda(s;z)$ series and using \eqref{ZetaLambdaComplement}, we obtain
		\begin{align}\label{ZetaLambdaContinuation}
			\begin{split}
				\zeta_\Lambda(s;{z})&=\sum_{\left|\lambda_k\right|\le \left|z_1\right|} \frac{1}{(\lambda_k-z_1)^s \dots (\lambda_k-z_n)^s}+\sum_{\left|\lambda_k\right|> \left|z_1\right|} \frac{1}{(\lambda_k-z_1)^s \dots (\lambda_k-z_n)^s}\\
				&=\sum_{\left|\lambda_k\right|\le \left|z_1\right|} \lambda_k^{-ns}\left(\prod_{j=1}^{n} \left(1-\frac{z_j}{\lambda_k}\right)^{-s}-\sum_{\ell=0}^{m} \frac{P_\ell(s;{z})}{\lambda_k^\ell}\right)
                \\&+\sum_{\left|\lambda_k\right|> \left|z_1\right|} \lambda_k^{-ns}\left(\prod_{j=1}^{n} \left(1-\frac{z_j}{\lambda_k}\right)^{-s}-\sum_{\ell=0}^{m} \frac{P_\ell(s;{z})}{\lambda_k^\ell}\right)\\
				&+\sum_{\ell=0}^{m}P_\ell(s;{z})\zeta_\Lambda(ns+\ell).
			\end{split}
		\end{align}
		The first sum in equation \eqref{ZetaLambdaContinuation} is finite and the second sum converges absolutely and uniformly on $\re(s)>\frac{\mu-m-1}{n}$. Therefore equation \eqref{ZetaLambdaContinuation} provides a meromorphic continuation of the function $\zeta_\Lambda(s;z)$ to $\re(s)>\frac{\mu-m-1}{n}$. This region contains $s=0$ as an holomorphic point. Thus proving that $(\Lambda-z_1)\dots(\Lambda-z_n)$ is regularizable. 
	\end{proof}

	\begin{theorem}\label{MainThm}
		The following equation holds
		\begin{align}\label{MainThmEq}
			\begin{split}
				D_\Lambda(z_1,\dots,z_n)&=\exp\left(-n\zeta_\Lambda'(0)-\sum_{j=1}^{n}\sum_{\ell=1}^{m} \frac{\fp \zeta_\Lambda(\ell)}{\ell} z_j^\ell-\frac{1}{n}\sum_{j=1}^{n} \sum_{\ell=2}^{m} \frac{\res \zeta_\Lambda(\ell) H_{\ell-1}}{\ell} z_j^\ell\right.\\
				&\qquad\quad\ \;-\frac{1}{n}\sum_{1\le j_1< j_2 \le n}\sum_{\ell=2}^{m}\sum_{r=1}^{\ell-1}\frac{\res\zeta_\Lambda(\ell)z_{j_1}^r z_{j_2}^{\ell-r}}{r(\ell-r)}\Bigg)\prod_{j=1}^{n} \Delta_\Lambda(z_j)\\
				&=\prod_{j=1}^{n} D_\Lambda(z_j)\exp\left(\left(1-\frac{1}{n}\right)\sum_{j=1}^{n} \sum_{\ell=2}^{m} \frac{\res \zeta_\Lambda(\ell) H_{\ell-1}}{\ell} z_j^\ell\right.\\
                &\qquad\qquad\qquad\qquad \ \left.-\frac{1}{n}\sum_{1\le j_1< j_2\le n}\sum_{\ell=2}^{m}\sum_{r=1}^{\ell-1}\frac{\res\zeta_\Lambda(\ell)z_{j_1}^r z_{j_2}^{\ell-r}}{r(\ell-r)}\right).
			\end{split}
		\end{align}
	\end{theorem}
	\begin{proof}
		Differentiating equation \eqref{ZetaLambdaContinuation} with respect to $s$ yields
		\begin{align*}
			\partial_s\zeta_\Lambda(s;{z})&=\sum_k \left(-\sum_{j=1}^{n}\log (\lambda_k-z_j)\right)\prod_{j=1}^{n}(\lambda_k-z_j)^{-s}-\sum_{\ell=0}^{m}\frac{\partial_s P_\ell(s;{z})-n\log \lambda_k P_\ell(s;{z})}{\lambda_k^{ns+\ell}}\\
			&+\sum_{\ell=0}^{m}\frac{d}{ds}P_\ell(s;{z})\zeta_\Lambda(ns+\ell)
		\end{align*}
		Substituting $s=0$ in above equation and using Lemma \ref{P_lLemma}, we obtain
		\begin{align}\label{OldZeta'0Delta}
			-\partial_s\zeta_\Lambda(0;{z})=\sum_k \left(\sum_{j=1}^{n}\log(\lambda_k-z_j)-n\log\lambda_k+\sum_{j=1}^{n}\sum_{\ell=1}^{m}\frac{z_j^\ell}{\ell \lambda_k^\ell} \right)-\sum_{\ell=0}^{m}\frac{d}{ds}P_\ell(s;{z})\zeta_\Lambda(ns+\ell)\Big\rvert_{s=0}.
		\end{align}
		Upon taking logarithms in equation \eqref{WeierstrassDef} and summing over $z_j$ for $j=1,\dots,n$, we get 
		\begin{align*}
			\sum_k \left(\sum_{j=1}^{n}\log(\lambda_k-z_j)-n\log\lambda_k+\sum_{j=1}^{n}\sum_{\ell=1}^{m}\frac{z_j^\ell}{\ell \lambda_k^\ell} \right)=\sum_{j=1}^{n}\log \Delta_\Lambda(z_j).
		\end{align*}
		Since $\zeta_\Lambda(s)$ has poles of order at most $1$, we have $\zeta_\Lambda(ns+\ell)=\frac{\res\zeta_\Lambda(\ell)}{ns}+\fp\zeta_\Lambda(\ell)+O(s)$ near $s=0$. $\frac{d}{ds}P_\ell(s;{z})\zeta_\Lambda(ns+\ell)\Big\rvert_{s=0}$ for $\ell=0,1$ are easily computed to be $n\zeta_\Lambda'(0)$ and $\left(\sum_{j=1}^n z_j\right)\fp\zeta_\Lambda(1)$ respectively. For $\ell\ge 2$, by Lemma \ref{P_lLemma}, we have
		\begin{align*}
			\frac{d}{ds}P_\ell(s;{z})\zeta_\Lambda(ns+\ell)\Big\rvert_{s=0}&=\frac{d}{ds}\left( p_{1,\ell}(z)s+p_{2,\ell}(z)s^2+O\left(s^3\right) \right)\left( \frac{\res\zeta_\Lambda(\ell)}{ns}+\fp\zeta_\Lambda(\ell)+O(s) \right)\Bigg\rvert_{s=0}\\
			&=p_{1,\ell}(z)\fp\zeta_\Lambda(\ell)+\frac{p_{2,\ell}(z)\res\zeta_\Lambda(\ell)}{n}.
		\end{align*}
		Therefore equation \eqref{OldZeta'0Delta} becomes
		\begin{align}\label{Zeta'0Delta}
			\begin{split}
				-\partial_s\zeta_\Lambda(0;{z})&=\sum_{j=1}^{n}\log \Delta_\Lambda(z_j)-n\zeta_\Lambda'(0)-\sum_{\ell=1}^{m}p_{1,\ell}(z)\fp \zeta_\Lambda(\ell)-\sum_{\ell=2}^{m}\frac{p_{2,\ell}(z)\res\zeta_\Lambda(\ell)}{n}\\
				&=\sum_{j=1}^{n}\log \Delta_\Lambda(z_j)-n\zeta_\Lambda'(0)-\sum_{j=1}^{n}\sum_{\ell=1}^{m} \frac{\fp \zeta_\Lambda(\ell)}{\ell} z_j^\ell-\frac{1}{n}\sum_{j=1}^{n} \sum_{\ell=2}^{m} \frac{\res \zeta_\Lambda(\ell) H_{\ell-1}}{\ell} z_j^\ell\\
				&-\frac{1}{n}\sum_{1\le j_1< j_2 \le n}\sum_{\ell=2}^{m}\sum_{r=1}^{\ell-1}\frac{\res\zeta_\Lambda(\ell)z_{j_1}^r z_{j_2}^{\ell-r}}{r(\ell-r)}
			\end{split}
		\end{align}
		Exponentiating the equation \eqref{Zeta'0Delta} and using the basic fact that $D_\Lambda(z_1,\dots,z_n)=\exp(-\partial_s\zeta_\Lambda(0;{z}))$ concludes the first part of the equation \eqref{MainThmEq}. By the fundamental theorem of zeta regularization (Theorem \ref{FundamentalThm}), we have
		\begin{align*}
			\log \Delta_\Lambda(z_j)-\zeta_\Lambda'(0)-\sum_{\ell=1}^{m}\frac{\fp \zeta_\Lambda(\ell)}{\ell} z_j^\ell=\log D_\Lambda(z_j)+\sum_{\ell=2}^{m} \frac{\res \zeta_\Lambda(\ell) H_{\ell-1}}{\ell} z_j^\ell
		\end{align*}
		for every $j=1,\dots,n$. Thus equation \eqref{Zeta'0Delta} becomes
		\begin{align}\label{Zeta'0RegProd}
            \begin{split}
                -\partial_s\zeta_\Lambda(0;{z})&=\sum_{j=1}^{n} \log D_\Lambda(z_j)+\left(1-\frac{1}{n}\right)\sum_{j=1}^{n} \sum_{\ell=2}^{m} \frac{\res \zeta_\Lambda(\ell) H_{\ell-1}}{\ell} z_j^\ell \\
                &-\frac{1}{n}\sum_{1\le j_1< j_2\le n}\sum_{\ell=2}^{m}\sum_{r=1}^{\ell-1}\frac{\res\zeta_\Lambda(\ell)z_{j_1}^r z_{j_2}^{\ell-r}}{r(\ell-r)}
            \end{split}
		\end{align}
		Similarly, exponentiating the equation \eqref{Zeta'0RegProd} concludes the second part of the equation \eqref{MainThmEq}. This completes the proof.
	\end{proof}
	
	\begin{corollary}
		The discrepancy $F=F(\Lambda-z_1,\dots,\Lambda-z_n)$ between the sequences $\Lambda-z_1,\dots,\Lambda-z_n$ is given by
		\begin{align*}
			F=\log \frac{\RegularProd_k \prod_{j=1}^{n} (\lambda_k-z_j)}{\prod_{i=j}^{n}\RegularProd_k (\lambda_k-z_j)}&=\left(1-\frac{1}{n}\right)\sum_{j=1}^{n} \sum_{\ell=2}^{m} \frac{\res \zeta_\Lambda(\ell) H_{\ell-1}}{\ell} z_j^\ell \\
            &-\frac{1}{n}\sum_{1\le j_1< j_2 \le n}\sum_{\ell=2}^{m}\sum_{r=1}^{\ell-1}\frac{\res\zeta_\Lambda(\ell)z_{j_1}^r z_{j_2}^{\ell-r}}{r(\ell-r)}.
		\end{align*}
	\end{corollary}
	\begin{corollary}
		If the sequence $\Lambda$ has order $\mu<2$, then the discrepancy $F$ vanishes and the following equation holds
		\begin{align*}
			\RegularProd_k \prod_{j=1}^{n} (\lambda_k-z_j)=\prod_{j=1}^{n}\RegularProd_k (\lambda_k-z_j)=\prod_{j=1}^{n} D_\Lambda(z_j).
		\end{align*}
	\end{corollary}
	\begin{corollary}
		Let $\varphi_j$ be given monic polynomials that don't vanish in $\Lambda$ for $j=1,\dots,n$ and $\Omega$ be the set of all roots of $\varphi_j$ counted with multiplicity. If the sequence $\Lambda$ has order $\mu<2$, then the following equation holds
		\begin{align*}
			\RegularProd_k \prod_{j=1}^{n} \varphi_j(\lambda_k)=\prod_{j=1}^{n}\RegularProd_k \varphi_j(\lambda_k)=\prod_{w\in\Omega}D_\Lambda(w).
		\end{align*}
	\end{corollary}
	By Theorem \ref{MainThm}, we see that the function $D_\Lambda(z_1,\dots,z_n)$ is an holomorphic function of $n$-variables in the whole complex plane. It is also worthwhile to note that while considering the polynomials $\lambda_k^n\pm w^n$, the sums $\sum_{j=1}^{n}z_j^\ell$ are easy to evaluate and often vanishing.
	\section{Applications And Examples}
	
	\subsection{Natural Numbers, Lerch-Wakayama-Mizuno Formulas}
    
	Let $\lambda_k=k$ for $k\ge 1$ and $\varphi_j$ be given monic polynomials that don't vanish in $\mathbb{Z}^+$ for $j=1,\dots,n$. Then the sequence $\Lambda=\mathbb{Z}^+$ has order $1$ and by Theorem \ref{MainThm}, the discrepancy $F(\mathbb{Z}^+-z_1,\dots,\mathbb{Z}^+-z_n)$ vanishes and so does $F(\varphi_1(\mathbb{Z}^+),\dots,\varphi_n(\mathbb{Z}^+))$. We therefore have the classical result of Mizuno.
	
	\begin{corollary}
		The following equation holds
		\begin{align*}
			\RegularProd_{k=1}^{\infty} \prod_{j=1}^{n} (k-z_j)=\prod_{j=1}^{n}\RegularProd_{k=1}^{\infty} (k-z_j)=\prod_{j=1}^{n}\frac{\sqrt{2\pi}}{\Gamma(1-z_j)}.
		\end{align*}
	\end{corollary}
	\begin{corollary}
		Let $\varphi_j$ be given monic polynomials that don't vanish in $\mathbb{Z}^+$ for $j=1,\dots,n$ and $\Omega$ be the set of all roots of $\varphi_j$ counted with multiplicity. Then the following equation holds
		\begin{align*}
			\RegularProd_{k=1}^{\infty} \prod_{j=1}^{n} \varphi_j(k)=\prod_{j=1}^{n}\RegularProd_{k=1}^{\infty} \varphi_j(k)=\prod_{w\in \Omega}\frac{\sqrt{2\pi}}{\Gamma(1-w)}.
		\end{align*}
	\end{corollary}
	
	\subsection{Multi-dimensional Products}\label{SectionMultidim}
    
	For a given positive integer $N$ and a vector $\omega=(\omega_1,\dots,\omega_N)\in\mathbb{C}^n$ where $\mathfrak{Re}(\omega_j)>0$ for $j=1,\dots,N$, we define the $N$th Barnes multiple zeta function by
	\begin{align*}
		\zeta_N(s,z;\omega)=\sum_{k\in\mathbb{N}^N}\frac{1}{(k\cdot \omega+z)^s}
	\end{align*}
	where $\mathfrak{Re}(s)>N$ and $\mathfrak{Re}(z)>0$. Thus we have $\mu=m=N$. Barnes \cite{Barnes1,Barnes2,Barnes3,Barnes4} has proved that the function $\zeta_N$ has a meromorphic continuation to the whole complex plane with only simple poles at $s=1,\dots,N$. The Barnes multiple Gamma function is defined as
	\begin{align*}
	    \Gamma_N(z;\omega)=\exp\left(\partial_s\zeta_N(0,z;\omega)\right)=\left(  \RegularProd_{k\in\mathbb{N}^N}\left( k\cdot \omega+z \right)\right)^{-1}.
	\end{align*}
        We remark that Barnes originally defined $\Gamma_N^B(z;\omega)=\rho_N(\omega)\Gamma_N(z;\omega)$ where $\rho_N(\omega)$ is the Barnes-Stirling modular gamma constant. The multiple gamma functions satisfy the fundamental functional equation
        \begin{align*}
            \frac{\Gamma_N(z+\omega_j;\omega)}{\Gamma_N(z;\omega)}=\Gamma_{N-1}(z;\omega(j))
        \end{align*}
        for every $j=1,\dots,N$ where $\omega(j)=\left( \omega_1,\dots,\omega_{j-1},\omega_{j+1},\dots,\omega_n \right)$. In the special case $\omega=(1,\dots,1)$, we obtain the classical Vigneras-Barnes mutliple zeta and gamma functions. Let multiple Bernoulli polynomials be given by
        \begin{align*}
            \frac{t^N e^{xt}}{\left( e^{\omega_1t}-1 \right)\dots\left( e^{\omega_Nt}-1 \right)}=\sum_{n=0}^{\infty}\frac{t^n}{n!}B_{N,n}(x;\omega)
        \end{align*}
        in a neighbourhood of $t=0$. It has been proven in \cite{OnBarnes} that
        \begin{align*}
            \res \zeta_N(\ell,z;\omega)=\frac{(-1)^{N-\ell}}{(\ell-1)!(N-\ell)!}B_{N,N-\ell}(z;\omega)
        \end{align*}
        where $\ell=1,\dots,N$. Let us denote
        \begin{align*}
            D_N(z_1,\dots,z_n;\omega)=\RegularProd_{k\in\mathbb{N}^N} \prod_{j=1}^{n} \left( k\cdot \omega-z_j \right)
        \end{align*}
        and
        \begin{align*}
            D_N(z;\omega)=\RegularProd_{k\in\mathbb{N}^N} \left( k\cdot \omega-z \right)=\frac{1}{\Gamma_N(-z;\omega)}.
        \end{align*}
        By the virtue of Theorem \ref{MainThm}, we have the following formula.
        \begin{corollary}
            The following equation holds
            \begin{align*}
                D_N(z_1,\dots,z_n;\omega)\prod_{j=1}^{n} \Gamma_N(-z_j;\omega)=\exp&\left(\left(1-\frac{1}{n}\right)\sum_{j=1}^{n} \sum_{\ell=2}^{N} \frac{(-1)^{N-\ell} H_{\ell-1}}{\ell!(N-\ell)!}B_{N,N-\ell}(0;\omega) z_j^\ell\right.\\
                & \left.-\frac{1}{n}\sum_{1\le j_1< j_2\le n}\sum_{\ell=2}^{N}\sum_{r=1}^{\ell-1}\frac{(-1)^{N-\ell}B_{N,N-\ell}(0;\omega)z_{j_1}^r z_{j_2}^{\ell-r}}{r(\ell-r)(\ell-1)!(N-\ell)!}\right).
            \end{align*}
        \end{corollary}
        If we change the variables $z_j\to -z_j$, we get
        \begin{align*}
            \RegularProd_{k\in\mathbb{N}^N} \prod_{j=1}^{n} \left( k\cdot \omega+z_j \right)\prod_{j=1}^{n} \Gamma_N(z_j;\omega)&=\exp\left(\left(1-\frac{1}{n}\right)\sum_{j=1}^{n} \sum_{\ell=2}^{N} \frac{(-1)^{N} H_{\ell-1}}{\ell!(N-\ell)!}B_{N,N-\ell}(0;\omega) z_j^\ell\right.\\
            &\qquad\qquad\left.-\frac{1}{n}\sum_{1\le j_1< j_2\le n}\sum_{\ell=2}^{N}\sum_{r=1}^{\ell-1}\frac{(-1)^{N}B_{N,N-\ell}(0;\omega)z_{j_1}^r z_{j_2}^{\ell-r}}{r(\ell-r)(\ell-1)!(N-\ell)!}\right)
        \end{align*}
        We give special attention to the case $N=2$. Mizuno \cite{Mizuno} has evaluated regularized products of form
        \begin{align*}
            \RegularProd_{k_1,k_2\in\mathbb{N}} \prod_{j=1}^{n} \left( k_1\omega_{1,j}+k_2\omega_{2,j}+z_j \right)
        \end{align*}
        where $\mathfrak{Re}(\omega_{1,j}),\mathfrak{Re}(\omega_{2,j}),\mathfrak{Re}(z_j)>0$ for every $j=1,\dots n$ with $\omega_{1,j}\neq\omega_{1,j'},\omega_{2,j}\neq\omega_{2,j'}$ and $\omega_{1,j}\omega_{2,j'}\neq\omega_{1,j'}\omega_{2,j}$ for $j\neq j'$. It is interesting to compare the results since in our case we have $\omega_{1,j}=\omega_{1,j'}=\omega_1$ and $\omega_{2,j}=\omega_{2,j'}=\omega_2$ for all $j,j'$. Mizuno has proved that
        \begin{align*}
            \RegularProd_{k_1,k_2\in\mathbb{N}} \prod_{j=1}^{n} \left( k_1\omega_{1,j}+k_2\omega_{2,j}+z_j \right)=e^{F_2}\prod_{j=1}^{n} \frac{1}{\Gamma_2(z_j;\omega)}
        \end{align*}
        where $F_2$ is given by
        \begin{align*}
            F_2=\frac{1}{2n}\sum_{1\le j_1< j_2 \le n}&\frac{\omega_{1,j_1}\omega_{2,j_2}-\omega_{1,j_2}\omega_{2,j_1}}{\omega_{1,j_1}\omega_{1,j_2}}\log\frac{\omega_{1,j_2}}{\omega_{1,j_1}}B_2\left(  \frac{\omega_{1,j_1}z_{j_2}-\omega_{1,j_2}z_{j_1}}{\omega_{1,j_1}\omega_{2,j_2}-\omega_{1,j_2}\omega_{2,j_1}}\right)\\
            &+\frac{\omega_{2,j_1}\omega_{1,j_2}-\omega_{1,j_1}\omega_{2,j_2}}{\omega_{2,j_1}\omega_{2,j_2}}\log\frac{\omega_{2,j_2}}{\omega_{2,j_1}}B_2\left(  \frac{\omega_{2,j_1}z_{j_2}-\omega_{2,j_2}z_{j_1}}{\omega_{2,j_1}\omega_{1,j_2}-\omega_{1,j_1}\omega_{2,j_2}}\right)
        \end{align*}
        and $B_2(x)=x^2-x+1/6$ is the second Bernoulli polynomial. Our results are indeed similar with Mizuno's. Taking $N=2$ in above equations, we obtain
        \begin{align*}
            \RegularProd_{k_1,k_2\in\mathbb{N}} \prod_{j=1}^{n} \left( k_1\omega_1+k_2\omega_2+z_j \right)\prod_{j=1}^{n} \Gamma_2(z_j;\omega)=\exp&\left(\left(1-\frac{1}{n}\right)\sum_{j=1}^{n} \frac{B_{2,0}(0;\omega)}{2}z_j^2\right.\\
            &\left.-\frac{1}{n}\sum_{1\le j_1< j_2\le n}B_{2,0}(0;\omega)z_{j_1}z_{j_2}\right).
        \end{align*}
        By a simple calculation, we have
        \begin{align*}
            B_{2,0}(0;\omega)=\lim_{t \to 0} \frac{t^2}{\left(e^{\omega_1t}-1\right)\left(e^{\omega_2t}-1\right)}=\frac{1}{\omega_1\omega_2}.
        \end{align*}
        Therefore finally, we get
        \begin{align*}
            \RegularProd_{k_1,k_2\in\mathbb{N}} \prod_{j=1}^{n} \left( k_1\omega_1+k_2\omega_2+z_j \right)\prod_{j=1}^{n} \Gamma_2(z_j;\omega)=\exp\left(\left(1-\frac{1}{n}\right)\sum_{j=1}^{n} \frac{z_j^2}{2\omega_1\omega_2}\right.\left.-\frac{1}{n}\sum_{1\le j_1< j_2\le n}\frac{z_{j_1}z_{j_2} }{\omega_1\omega_2}\right).
        \end{align*}

        \subsection{Zeros Of The Riemann Zeta Function}
        
        The Riemann Zeta function is defined as
        \begin{align*}
            \zeta(s)=\sum_{n=1}^{\infty}\frac{1}{n^s} \qquad (\mathfrak{Re}(s)>1).
        \end{align*}
        It can be analytic continued to the entire complex plane except $s=1$ where it has a simple pole. The Riemann Xi function is defined as
        \begin{align*}
            \xi(s)=\frac{1}{2}s(s-1)\pi^{-s/2}\Gamma\Big(\frac{s}{2}\Big)\zeta(s)
        \end{align*}
        and it satisfies the functional equation $\xi(s)=\xi(1-s)$ \cite{TheoryofZeta}. The function $\xi(s)$ is an entire function of order $1$ and has the well-known Weierstrass product form 
        \begin{align*}
            \xi(s)=\frac{1}{2}\prod_{\rho}\Big(1-\frac{s}{\rho}\Big)
        \end{align*}
        where the product is taken over all zeros of the function $\xi$, which are non-trivial zeros of the Riemann Zeta function. The zeta function associated to the non-trivial zeros of the Riemann zeta function is defined as
        \begin{align*}
            Z(s)=\sum_{\rho}\frac{1}{\rho^s}, \qquad (\mathfrak{Re}(s)>1)
        \end{align*}
        and has been thoroughly studied in \cite{VorosBook}. The function $Z$ has the particular values $Z(0)=2$
        $Z'(0)=\frac{1}{2}\log 2$ and $\mathcal{FP}Z(1)=1-\frac{1}{2}\log 2+\frac{1}{2}C$ where $C$ is the Euler's constant. By Theorem \ref{FundamentalThm}, we have
        \begin{align*}
            \RegularProd_{\rho}(\rho-s)=\sqrt{2}\xi(s)   \qquad \text{or}\qquad  \RegularProd_{\rho}(\rho+s)=\sqrt{2}\xi(s+1).
        \end{align*}
        Since the non-trivial zeros of the zeta function form a sequence of order $1$, the discrepancy between polynomials vanishes identically.
        \begin{corollary}
            The following equation holds
            \begin{align*}
                \RegularProd_{\rho}\prod_{j=1}^{n}(\rho+z_j)=\prod_{j=1}^{n}\sqrt{2}\xi(1+z_j)=\prod_{j=1}^{n}\RegularProd_{\rho}(\rho+z_j).
            \end{align*}
        \end{corollary}
        \begin{corollary}
		Let $\varphi_j$ be given monic polynomials that don't vanish at the non-trivial zeros of the zeta function for $j=1,\dots,n$ and $\Omega$ be the set of all roots of $\varphi_j$ counted with multiplicity. Then, the following equation holds
		\begin{align*}
			\RegularProd_\rho \prod_{j=1}^{n} \varphi_j(\rho)=\prod_{j=1}^{n}\RegularProd_\rho \varphi_j(\rho)=\prod_{w\in\Omega}\sqrt{2}\xi(w).
		\end{align*}
	\end{corollary}
    Let $\gamma=\gamma(\rho)$ be determined as $\rho=\frac{1}{2}+i\gamma$. Using the property $\RegularProd a\lambda_k=a^{\zeta_\Lambda(0)}\RegularProd \lambda_k$ and $Z(0)=2$, we obtain
    \begin{align*}
        \RegularProd_{\gamma}(\gamma-s)=-\sqrt{2}\xi\Big(\frac{1}{2}+is\Big)\qquad \text{or}\qquad \RegularProd_{\gamma}(\gamma+s)=-\sqrt{2}\xi\Big(\frac{1}{2}-is\Big).
    \end{align*}
    Similarly, the discrepancy vanishes and we get the following.
    \begin{corollary}
        The following equation holds
        \begin{align*}
            \RegularProd_{\gamma}\prod_{j=1}^{n}(\gamma+z_j)=\prod_{j=1}^{n}-\sqrt{2}\xi\left(\frac{1}{2}-iz_j\right)=\prod_{j=1}^{n}\RegularProd_{\rho}(\gamma+z_j).
        \end{align*}
    \end{corollary}
    \begin{corollary}
		Let $\varphi_j$ be given monic polynomials that don't vanish at $\gamma$ for $j=1,\dots,n$ and $\Omega$ be the set of all roots of $\varphi_j$ counted with multiplicity. Then, the following equation holds
		\begin{align*}
			\RegularProd_\gamma \prod_{j=1}^{n} \varphi_j(\gamma)=\prod_{j=1}^{n}\RegularProd_\gamma \varphi_j(\gamma)=\prod_{w\in\Omega}-\sqrt{2}\xi\Big(\frac{1}{2}+iw\Big).
		\end{align*}
	\end{corollary}
    
    \subsection{Zeros Of The Bessel Function}
    
    Let $\nu\ge \frac{1}{2}$ be a real number. The Bessel function of the first kind $J_\nu$ is given by
    \begin{align*}
        J_\nu(z)=\sum_{n=0}^{\infty}\frac{(-1)^n}{n!\Gamma(n+\nu+1)}\left( \frac{z}{2} \right)^{2n+\nu}.
    \end{align*}
    Lommel's theorems (see \cite{Watson}) state that the function $J_\nu$ for $\nu>-1$ has only real zeros and furthermore there are an infinite amount of real zeros. Let us denote the positive zeros of the function $J_\nu$ by $0<j_{\nu,1}<\dots<j_{\nu,k}<\dots$. The zeros together with their negatives $-j_{\nu,k}$ form a sequence of order $1$. The Bessel zeta function, which was studied systematically for the first time by Hawkins in \cite{Hawkins}, is defined as
    \begin{align*}
        \zeta_{J_\nu}(s)=\sum_{k=1}^{\infty}\left( \frac{j_{\nu,k}}{\pi} \right)^{-s}, \qquad (\mathfrak{Re}(s)>1).
    \end{align*}
    For an extensive study of the function $\zeta_{J_\nu}$ we refer the reader to \cite{Actor} and references therein. We consider products over all zeros of the Bessel function $J_\nu$ except $0$. Thus we study the associated zeta function
    \begin{align*}
        \zeta^*_{J_\nu}(s)=\left( 1+e^{-\pi is} \right)\zeta_{J_\nu}(s).
    \end{align*}
    By the classical product formula for the Bessel function, we have the Weierstrass product
    \begin{align*}
        \Delta^*_{J_\nu}(z)=\prod_{k=1}^{\infty}\left( 1-\frac{\pi^2 z^2}{j_{\nu,k}^2} \right)=\left( \frac{\pi z}{2} \right)^{-\nu}\Gamma(\nu+1)J_\nu(\pi z).
    \end{align*}
    It has been proved in \cite{Actor} that $\zeta_{J_\nu}(0)=-\frac{1}{2}\left(\nu+\frac{1}{2}\right)$ and
    \begin{align*}
        \zeta_{J_\nu}'(0)=\frac{1}{2}\log\frac{2^{\nu-\frac{1}{2}}\Gamma(\nu+1)}{\pi^{\nu+1}}.
    \end{align*}
    Therefore we have
    \begin{align*}
        \partial_s\zeta^*_{J_\nu}  (0)=2\zeta_{J_\nu}'(0)-\pi i\zeta_{J_\nu}(0)=\log\frac{2^{\nu-\frac{1}{2}}\Gamma(\nu+1)}{\pi^{\nu+1}}+\frac{\pi i}{2}\left(\nu+\frac{1}{2}\right).
    \end{align*}
    The quantity $\mathcal{FP}\zeta_{J_\nu}(1)$ does not seem to have a closed form. By Theorem \ref{FundamentalThm}, we obtain
    \begin{align*}
        D_{J_\nu}^*(z)=\RegularProd_{k\neq0}\left( \frac{j_{\nu,k}}{\pi}-z \right)&=\Delta_{J_\nu}^*(z)\exp\left( -\partial_s\zeta^*_{J_\nu}(0)-z\mathcal{FP}\zeta^*_{J_\nu}(1) \right)\\
        &=\frac{\sqrt{2} \pi}{z^\nu}J_{\nu}(\pi z)\exp\left( -\frac{\pi i}{2}\left(\nu+\frac{1}{2}\right)-z\mathcal{FP}\zeta^*_{J_\nu}(1) \right)
    \end{align*}
    One might also use the property $\RegularProd a\lambda_k=a^{\zeta_\Lambda(0)}\RegularProd \lambda_k$ and $\zeta^*_{J_\nu}(0)=2\zeta_{J_\nu}(0)=-\left(\nu+\frac{1}{2}\right)$ to get
    \begin{align*}
        \RegularProd_{k\neq0}(j_{\nu,k}-z)=\frac{\sqrt{2\pi} }{z^\nu}J_{\nu} (z)\exp\left( -\frac{\pi i}{2}\left(\nu+\frac{1}{2}\right)-\frac{z}{\pi}\mathcal{FP}\zeta^*_{J_\nu}(1) \right).
    \end{align*}
    Taking $z=0$, we have the special cases
    \begin{gather*}
        \RegularProd_{k\neq0}\frac{j_{\nu,k}}{\pi}=\frac{\pi^{\nu+1} }{2^{\nu-\frac{1}{2}}\Gamma(\nu+1)}\exp\left( -\frac{\pi i}{2}\left(\nu+\frac{1}{2}\right) \right)\\
        \RegularProd_{k\neq0}j_{\nu,k}=\frac{\sqrt{2\pi} }{2^{\nu}\Gamma(\nu+1)}\exp\left( -\frac{\pi i}{2}\left(\nu+\frac{1}{2}\right) \right).
    \end{gather*}
    None of the formulas above seem to be recorded in the literature before. The zeros form a sequence of order $1$ and thus the discrepancy vanishes identically.
    \begin{corollary}
        The following equation holds
        \begin{align*}
            \RegularProd_{k\neq0}\prod_{j=1}^{n}\left( \frac{j_{\nu,k}}{\pi}-z_j \right)&=\prod_{j=1}^{n}\RegularProd_{k\neq0}\left( \frac{j_{\nu,k}}{\pi}-z_j \right)\\
            &=\prod_{j=1}^{n}\frac{\sqrt{2} \pi}{z_j^\nu}J_{\nu}(\pi z_j)\exp\left( -\frac{\pi i}{2}\left(\nu+\frac{1}{2}\right)-z_j\mathcal{FP}\zeta^*_{J_\nu}(1) \right).
        \end{align*}    
    \end{corollary}
    \begin{corollary}
        Let $\varphi_j$ be given monic polynomials that don't vanish at $j_{\nu,k}/\pi$ for $j=1,\dots,n$ and $\Omega$ be the set of all roots of $\varphi_j$ counted with multiplicity. Then, the following equation holds
        \begin{align*}
            \RegularProd_{k\neq0}\prod_{j=1}^{n}\varphi_j\left( \frac{j_{\nu,k}}{\pi} \right)&=\prod_{j=1}^{n}\RegularProd_{k\neq0}\varphi_j\left( \frac{j_{\nu,k}}{\pi} \right)\\
            &=\prod_{w\in\Omega}\frac{\sqrt{2} \pi}{w^\nu}J_{\nu}(\pi w)\exp\left( -\frac{\pi i}{2}\left(\nu+\frac{1}{2}\right)-w\mathcal{FP}\zeta^*_{J_\nu}(1) \right).
        \end{align*}    
    \end{corollary}
    We note that many of the special functions satisfying second order ordinary linear differential equations have analogous properties discussed here. Therefore this process is not special to the Bessel functions and may be done for a variety of special functions. One example that is worthwhile to note is the Airy function. 


\end{document}